\documentclass[twoside,a4,10pt]{amsart}
\usepackage{amssymb,amsmath,graphicx}

\def\Res{\mathop{\rm res}}

\def\R{\mathbb{R}}

\def\C{\mathbb{C}}
\def\CC{\mathfrak{C}}

\def\Re{{\rm Re} \,}

\def\2int{\mathop{\int\int}}

\def\kzu1#1{\buildrel #1 \over \longrightarrow}

\def\ds{\displaystyle}
\def\ts{\textstyle}

\def\diag{{\rm diag\,}}

\def\und{\quad{\rm and}\quad}

\def\v{\mathfrak{v}}
\def\V{\mathfrak{V}}

\def\0{\mathbf{0}}
\def\1{\mathbf{1}}

\def\be#1{\begin{equation}\label{#1}}
\def\ee{\end{equation}}
\def\ende{{\square}}
\def\Ende{\hfill$\ende$}

\def\wcirc{\buildrel \circ\over w}
\begin{document}

\title{Laplace contour integrals and linear differential equations}
\author{Norbert Steinmetz}
\maketitle

{\small\begin{center}{\it Dedicated to the memory of Professor Walter Hayman}\bigskip

{\sc Abstract}\bigskip

\begin{minipage}{100mm}The purpose of this paper is to determine the main properties of Laplace contour integrals
$$\Lambda(z)=\frac1{2\pi i}\int_\CC\phi_L(t)e^{-zt}\,dt,$$
that solve linear differential equations
$$L[w](z):=w^{(n)}+\sum_{j=0}^{n-1}(a_j+b_jz)w^{(j)}=0.$$
This concerns, in particular, the order of growth, asymptotic expansions, the Phragm\'en-Lindel\"of indicator, the  distribution of zeros,
the existence of sub-normal and polynomial solutions, and the corresponding Nevanlinna functions.\bigskip


\end{minipage}\end{center}}

\section{Introduction}
 Special functions usually admit several quite different definitions and representations. For example,
{\it Airy's equation}
$$w''-zw=0$$
has a particular solution, known as {\it Airy function}, that may be written as {\it Laplace contour integral}
\be{AiryLaplace}{\rm Ai}(z)=\frac1{2\pi i}\int_\CC e^{-zt+t^3/3}\,dt;\ee
the contour $\CC$ consists of the straight line from $+\infty\,e^{i\pi/3}$ to the origin followed by the straight line from $0$ to $+\infty\,e^{-i\pi/3}$.
In \cite{GHW} the authors G.~Gundersen, J.~Heittokangas and Z-T.~Wen  investigated two special families of linear differential equations,
namely the three-term equations
\be{FAM1}L[w]:=w^{(n)}+(-1)^{n+1}bw^{(k)}+(-1)^{n+1}zw=0\quad(0<k<n,~b\in\C)\ee
and
\be{FAM2}L[w]:=w^{(n)}-zw^{(k)}-w=0\quad(1<k<n)\ee
with the intention
{\it ``to make a contribution to this topic which includes a generalization of the Airy integral $[\ldots]$ to have more examples
of solutions of complex differential equations that have concrete properties''}.
To this end they determined several contour integral solutions
$$\frac1{2\pi i}\int_\CC F(z,w)\,dw$$
with varying kernels $F$ and contours $\CC$. Although \cite{GHW} does not contain any hint how to find appropriate kernels and contours,
not to mention appropriate (families of) differential equations like (\ref{FAM1}) and (\ref{FAM2}), the examples in \cite{GHW} reveal the nature of the contour integrals: in any case they are equal or may be transformed into {\it Laplace contour integrals}
\be{LaplaceTransform}\frac1{2\pi i}\int_\CC \phi(t)e^{-zt}\,dt\ee
of particular analytic functions $\phi$ over canonical contours or paths of integration $\CC$.
We will prove that to each linear differential equations
\be{GenLinCase}L[w]:=w^{(n)}+\sum_{j=0}^{n-1}(a_j+b_jz)w^{(j)}=0\quad(a_0+b_0z\not\equiv 0)\ee
there exists some distinguished Laplace contour integral solution (\ref{LaplaceTransform})
with kernel $\phi_L$ that is uniquely determined by the operator $L$ and itself determines canonical contours $\CC$. The main properties of these solutions -- denoted $\Lambda_L$ --,
which strongly remind on the Airy integral,
will be revealed. This concerns, in particular, the order of growth, the deficiency of the value zero, asymptotic expansions in particular sectors, the distribution of zeros, the Phragm\'en-Lindel\"of indicator, the Nevanlinna functions $T(r,\Lambda_L)$ and $N(r,1/\Lambda_L)$, and the existence of sub-normal solutions.

\section{Kernels and contours}

The following reflections on linear differential equations
\be{Linearoperator}L[w]:=\sum_{j=0}^nP_j(z)w^{(j)}=0,\ee
are more or less part of {\it mathematical folklore}, but unfortunately not as well-known as they should be. One can find few remarks in Wasow~\cite{Wasow}, p.~123~ff.,
where the method of Laplace contour integrals and the {\it saddle-point method} is applied to Airy's equation. Hille~\cite{Hille}, p.~216~ff., and Ince~\cite{ince}, chapter XVIII, deal with linear differential equations (\ref{Linearoperator})
under the constraint $\deg P_n\ge \deg P_j$. This, however, is by no means necessary and even prevents the discussion of  the most important case $P_n(z)\equiv 1$,
where the solutions are {\it entire functions of finite order of growth}.\smallskip

We start with linear differential equations (\ref{Linearoperator}) with polynomial coefficients
\be{Coefficients}P_j(z)=\sum_\alpha c_{j\alpha}z^\alpha\quad(P_n(z)P_0(z)\not\equiv 0,~d=\max\deg P_j>0)\ee
and are looking for solutions that may be written as Laplace contour integrals (\ref{LaplaceTransform}); the function $\phi$ and thereby possible paths of integration $\CC$ have to be determined.
Formally we have
$$w^{(j)}(z)=\frac{1}{2\pi i}\int_\CC \phi_j(t)e^{-zt}\,dt\quad{\rm with}~\phi_j(t)=(-t)^j\phi(t)$$
and
\be{Teil}\begin{aligned}
(-1)^\alpha\int_\CC z^\alpha \phi_j(t)e^{-zt}\,dt=&~\int_\CC \phi_j(t)D^\alpha [e^{-zt}]\,dt\qquad({\rm where~}D= d/{dt})\cr
=&~\big[ \phi_j(t)D^{\alpha-1} [e^{-zt}]\big]_\CC-\int_\CC \phi_j'(t)D^{\alpha-1}[e^{-zt}]\,dt\end{aligned}\ee
if $\alpha\ge 1$; here
$$\big[\Psi(t)\big]_\CC=\lim_{t\to b}\Psi(t)-\lim_{t\to a}\Psi(t)$$
denotes the variation of $\Psi$ along $\CC$, when $\CC$ starts at $a$ and ends at $b$; of course, $\CC$ may be closed. We demand
\be{Rand}\big[\phi^{(\alpha-1)}_j(t) e^{-zt}\big]_\CC=0\quad(0\le j\le n,~1\le\alpha\le\max\limits_{0\le j\le n}\deg P_j)\ee
to obtain
$$L[w](z)=\frac1{2\pi i}\int_\CC\Big[\sum_{j,\alpha}c_{j\alpha}D^\alpha[(-t)^j\phi(t)]\Big] e^{-zt}\,dt,$$
hence $L[w](z)\equiv 0$ if
\be{SuffCond}\sum_{j,\alpha}c_{j\alpha}D^\alpha[(-t)^j\phi(t)]\equiv 0.\ee
These calculations are justified if (\ref{Rand}) holds and the integrals converge absolutely with respect to arc-length, and locally uniformly with respect to $z$. With
\be{Qalpha}Q_\alpha(t)=\sum_{j=0}^nc_{j\alpha}(-t)^j,\ee
condition (\ref{SuffCond}) may be written as
\be{SuffCond2}\sum_{\alpha=0}^d \big[Q_\alpha\phi\big]^{(\alpha)}(t)\equiv 0\quad(d=\max\deg P_j),\ee
and we have the following theorem, which forms the basis of our considerations.

\subsection{\sc Theorem}\label{MainTheo} {\it Let $\phi$ be any non-trivial solution to $(\ref{SuffCond2})$ and $\CC$ any contour such that the integrals $(\ref{Teil})$
converge absolutely with respect to arc-length of $\CC$ and
locally uniformly with respect to $z$, and also $(\ref{Rand})$ holds. Then the Laplace contour integral $(\ref{LaplaceTransform})$
solves the differential equation $(\ref{Linearoperator})$.}\medskip

Of course one has to check in each particular case that $w$ is non-trivial. In most cases this is a corollary of the inverse formula and/or uniqueness theorems for the Laplace transformation.

\section{Equations with coefficients of degree one}

\subsection{\sc The kernel}  The method of the previous section turns out to be most useful in the study of equation
(\ref{GenLinCase}), in which case
$$Q_0(t)=(-t)^n+\sum_{j=0}^{n-1}a_j(-t)^j\und Q_1(t)=b_q(-t)^q+\sum_{j=0}^{q-1}b_j(-t)^j,$$
where $b_q\ne 0$ and $b_j=0$ for $j>q$ is assumed. By a linear change of the independent variable $z$ we may achieve
\be{bqNormalisation}b_q=(-1)^{n-q+1}\quad({\rm and~} a_q=0~{\rm if~desired});\ee
this will be assumed henceforth. Then $Q_0\phi+(Q_1\phi)'=0$ gives
\be{phiford1}\phi(t)=\frac1{Q_1(t)}\exp\Big[-\int\frac{Q_0(t)}{Q_1(t)}\,dt\Big]\ee
up to an arbitrary constant factor.  More precisely,
\be{Generalphi}
\begin{array}{c}\phi(t)=\ds\prod_{\nu=1}^m(t-t_\nu)^{-m_\nu-\lambda_\nu}\exp\Big[R_0(t)+\sum_{\nu=1}^m R_\nu\Big(\frac1{t-t_\nu}\Big)\Big]\quad{\rm with}\cr
\ds\lambda_\nu=\Res_{t_\nu}\Big[\frac{Q_0}{Q_1}\Big],\quad \deg R_\nu\le m_\nu-1, \quad R_0(t)=\frac{t^{n-q+1}}{n-q+1}+\cdots\end{array}\ee
holds. Sometimes it suffices to know that $\phi$ has the form
\be{phinormal}\phi(t)=\exp\Big[\frac{t^{n-q+1}}{n-q+1}+\psi(t)\Big],\ee
where $\psi$ is holomorphic on $\{t:|t|>R_0,~ |\arg t|<\pi\}$ and satisfies
\be{psinormal}|\psi(t)|\le C|t|^{n-q}\quad{\rm as}~t\to\infty.\ee
Occasionally, any function (\ref{phiford1}) will denoted $\phi_L$.

\subsection{\sc Appropriate contours} There is some freedom in the choice of the contour $\CC$,
which is restricted only by condition (\ref{Rand}). Closed simple contours work in particular cases only, namely when $\lambda_\nu$ is an integer, and yield less interesting examples. More interesting are the contours
$\CC_{R,\alpha,\beta}$ that consist of three arcs as follows:
\begin{itemize}
\item[1.] the line $t=re^{i\alpha}$, where $r$ runs from $+\infty$ to $R\ge 0$;
\item[2.] the circular arc on $|t|=R$  from $Re^{i\alpha}$ via $t=Re^{i(\alpha+\beta)/2}$ to $Re^{i\beta}$;
\item[3.] the line $t=re^{i\beta}$, where $r$ runs from $R$ to $+\infty$.
\end{itemize}
In place of $\CC_{0,\alpha,\beta}$\,, $\CC_{R,\alpha,-\alpha}$\,, and $\CC_{0,\alpha,-\alpha}$ we will also write $\CC_{\alpha,\beta}$\,, $\CC_{R,\alpha}$\,, and $\CC_{\alpha}$\,,  respectively.
Note that $\CC_{R,\alpha,-\alpha}$ passes through $R$, while $\CC_{R,\alpha,2\pi-\alpha}$ passes through $-R$.
To ensure convergence of the contour integral and condition (\ref{Rand}), the real part of $t^{n-q+1}$ has to be negative on $\arg t=\alpha$ and $\arg t=\beta$;
this will tacitely be assumed. A canonical and our preferred choice is $\alpha=\frac{\pi}{n-q+1}$ and $\beta=-\alpha$. It is, however, almost obvious by Cauchy's theorem that the contour integral is independent of $\alpha,$ $\beta,$ and $R$ within their natural limitations. This will be proved later on.

\subsection{\sc The protagonist} By $w=\Lambda_{L}(z)$ we will denote any solution to equation (\ref{GenLinCase}) given by
\be{LAMBDA}\Lambda_L(z)=\frac1{2\pi i}\int_{\CC_{R,\frac{\pi}{n-q+1}}}\phi_L(t)e^{-zt}\,dt\ee
with associated function (\ref{phiford1}). Of course, $\Lambda_L$ is (up to a constant factor) uniquely determined by the operator $L$ and differs from operator to operator, but its essential properties depend only on $n-q$. Hence the Airy function is a typical representant of the solutions $\Lambda_L$ with $q=n-2$, where $n$ may be arbitrarily large.

\subsection{\sc Elementary Examples}

\medskip\begin{itemize}
\item[1)] {Airy's Equation} $w''-zw=0$ and $\Lambda_L={\rm Ai}$:\medskip

\noindent $Q_0(t)=t^2,$ $Q_1(t)=-1$, $\phi_L(t)=e^{t^3/3}$, and $\CC=\CC_{\frac\pi 3}$.
By Cauchy's theorem, any contour $\CC_{R,\alpha,\beta}$ with $R\ge 0$, $\alpha\in (\frac\pi 6,\frac\pi 2)$, and $\beta\in (-\frac\pi 6,-\frac\pi 2)$,
is admissible.\medskip

\item[2)] $w^{(n)}+(-1)^{n+1}bw^{(k)}+(-1)^{n+1}zw=0$, see (\ref{FAM1}):\medskip

\noindent $Q_0(t)=(-1)^nt^n-(-1)^{n+1+k}bt^k$, $Q_1(t)=(-1)^n$, hence
$$\phi_L(t)=\exp\Big[\frac{t^{n+1}}{n+1}+\frac{(-1)^{k+1}bt^{k+1}}{k+1}\Big]$$
(see formula (3.2) in \cite{GHW}).The contour $\CC=\CC_{\frac\pi{n+1}}$ is canonical.\medskip

\item[3)] $w^{(n)}-zw^{(k)}-w=0$ ($1<k<n$), see (\ref{FAM2}): \medskip

\noindent $Q_0(t)=(-1)^nt^n-1$, $Q_1(t)=(-1)^{k+1}t^k$, hence
$$\phi_L(t)=\frac1{t^k}\exp\Big[\frac{(-1)^{n-k}t^{n-k+1}}{n-k+1}+\frac{(-1)^kt^{-k+1}}{k-1}\Big].$$
Since $\phi_L$ has an essential singularity at $t=0$, one may choose $\CC:|t|=1$; it is, however, more interesting to choose $\CC=\CC_{{\frac\pi{n-k+1}}}$
if $n-q$ is even; if $n-q$ is odd, the differential equation does match our normalisation only after the change of variable $z\mapsto -z$.
We note that the substitution (the conformal map) $t=-1/u$ transforms the Laplace contour integral with kernel $\phi_L$ up to sign into some integral
$$\frac1{2\pi i}\int_{\tilde\CC}u^{k-2}\exp\Big[\frac zu-\frac{u^{-n+k-1}}{n-k+1}-\frac{u^{k-1}}{k-1}\Big]\,du$$
in accordance with formula (5.2) in \cite{GHW}.\medskip
\end{itemize}

\section{Asymptotic expansions and the order of growth}

\subsection{\sc The characteristic equation}\label{CharEq} The information on possible orders of growth and asymptotic expansions of transcendental solutions to (\ref{GenLinCase}) is
encoded in the algebraic equation
\be{AlgEqu}y^n +\sum_{j=0}^{n-1}(a_j+b_jz)y^j=0,\ee
which justifiably may be called {\it characteristic equation}.
Remember that we assume $b_q=(-1)^{n-q+1}$ and $b_j=0$ for $j>q$. Also let $p\le q$ denote the smallest index such that $b_p\ne 0$.
As $z\to\infty$, equation (\ref{AlgEqu}) has solutions\medskip

\begin{itemize}
\item[1)] $y_j(z)\sim \gamma_j z^{1/(n-q)}$ ($\gamma_j^{n-q}=(-1)^{n-q}$, $1\le j\le n-q$), hence
$$\int y_j(z)\,dz\sim \frac{\gamma_jz^{1+1/(n-q)}}{1+1/(n-q)}.$$
in any case ($n-q$ odd or even), one $\gamma_j$ equals $-1$.

\item[2)] $y_{n-q+j}\sim t_j$ ($1\le j\le q-p)$, where $t_j\ne 0$ is a root of $Q_1(-t)=0$ (counting multiplicities), hence
$$\int y_{n-q+j}(z)\,dz\sim t_j z;$$

\item[3)] $y_{n-p+j}(z)\sim \tau_jz^{-1/p}$ ($\tau_j^p=-a_0/b_p$, $1\le j\le p)$, hence
$$\qquad\quad\int y_{n-p+j}(z)\,dz\sim \frac{\tau_jz^{1-1/p}}{1-1/p}~{\rm if}~ p>1 {\rm ~and~}\int y_n(z)\,dz\sim\tau_1\log z~{\rm if~} p=1.$$
\end{itemize}
By $z^{1/r}$ and $\log z$ we mean the branches on $|\arg z|<\pi$ that are real on the positive real axis. If $p=1$, polynomial solutions may exist, but not otherwise.

\subsection{\sc The order of growth} Let
$$f(z)=\sum_{n=0}^\infty c_nz^n$$
be any transcendental entire function with {\it maximum modulus}, {\it maximum term}, and {\it central index}
$$\begin{aligned}
M(r,f)=&~\max\{|f(z)|:|z|=r\},\cr
\mu(r)=&~\max\{|c_n|r^n:n\ge 0\},{\rm ~and~}\cr
\nu(r)=&~\max\{n:|c_n|r^n=\mu(r)\},\end{aligned}$$
respectively. Then $f$ has order of growth
$$\varrho(f)=\limsup_{r\to\infty}\frac{\log\log M(r,f)}{\log r}=\limsup_{r\to\infty}\frac{\log\log\mu(r)}{\log r}=\limsup_{r\to\infty}\frac{\log\nu(r)}{\log r}.$$
The {\it central index method} for solutions to linear differential equations (\ref{Linearoperator}) is based on the relation
$$\frac{f^{(j)}(z_r)}{f(z_r)}=\Big(\frac{\nu(r)}{z_r}\Big)^j(1+\epsilon_j(r))\quad(j=1,2,\ldots,),$$
where $z_r$ is any point on $|z|=r$ such that $|f(z_r)|=M(r,f)$. This holds with $\epsilon_j(r)\to 0$ as $r\to\infty$
with the possible exception of some set $E_j$ of finite logarithmic measure; for a proof see Wittich~\cite{Wittich}, chapter I, or Hayman~\cite{hayman2}.
If $w=f(z)$ solves (\ref{GenLinCase}), $\nu(r)/z_r=v(r)$ satisfies
$$(1+\epsilon_n(r)v(r)^n+\sum_{j=0}^{n-1}(a_j+b_jz_r)(1+\epsilon_j(r))v(r)^j=0,$$
hence any such $v(r)$ is asymptotically correlated with some solution to the characteristic equation (\ref{AlgEqu}) as follows:
$$\nu(r)\sim z_r y_j(z_r)$$
as $r\to\infty$ holds even without exceptional set by the regular behaviour of $z_ry_j(z_r)$; since $\nu(r)$ is positive,  $z_ry_j(z_r)$ is asymptotically positive,
this giving  additional information on the possible values of $\arg z_r$.
Thus the transcendental solutions to equation (\ref{GenLinCase}) have possible orders of growth\medskip

\begin{tabular}{rl}$\varrho=1+\frac1{n-q}$:& solutions of this order always exist;\cr
$\varrho=1$:& necessary  for solutions of this order to exist  is $p<q$;\cr
$\varrho=1-\frac1p$:& necessary for solutions this order to exist  is $p>1$;\cr
$\varrho=0$:& necessary for polynomial solutions is $p=1$.\end{tabular}

\subsection{\sc Asymptotic expansions}
The differential equation (\ref{GenLinCase}) may  be rewritten in the usual way as a first-order system
\be{System}z^{-1}\v'=\big(B+z^{-1}A\big)\v\ee
for $\v=(w,w',\ldots,w^{(n-1)})$ with $n\times n$-matrices $A$ and $B$.
The following details can be found in Wasow's fundamental monograph~\cite{Wasow}. The system (\ref{System}) has {\it rank} one, but the matrix $B$ has vanishing eigenvalues.
Thus the theory of asymptotic integration yields only a local result (Theorem~19.1. in~\cite{Wasow}), which makes its applicability unpleasant and involved: to every angle $\theta$ there exists some sector $S_\theta:|\arg(ze^{-i\theta})|<\delta$ such that the system (\ref{System}) has a distinguished {\it fundamental matrix}
$$\V(z)=V(z|\theta)z^Ge^{\Pi(z)};$$                                                                                                                                                                                                since $B+z^{-1}A$ is holomorphic on $|z|>0$, the number $\delta$, the constant $n\times n$-matrix $G$, and the matrix $\Pi=\diag(\Pi_1,\ldots,\Pi_n)$, whose entries are polynomials in $z^{1/r}$ for some positive integer $r$, are universal, that is, they do not depend on $\theta$; $V(z|\theta)$ has an asymptotic expansion
\be{Vasymp}V(z|\theta)\sim \sum_{j=0}^\infty V_j\,z^{-1/r}\quad{\rm on~} S_\theta;\ee
the latter means
$$\|V(z|\theta)-\sum_{j=0}^m V_j\,z^{-1/r}\|=o(|z|^{-m/r})\quad(m=0,1,2,\ldots)$$
as $z\to\infty$ on $S_\theta$, uniformly on every closed sub-sector; again the matrices $V_j$ are independent of $\theta$. The fundamental matrices $\V$, that is, the matrices $V(z|\theta)$ may vary from sector to sector.[\footnote{Independence of the matrices $V_j$ and dependence of $V(z|\theta)$ on $\theta$ yields no contradiction. Asymptotic series may represent different analytic functions on sectors.}]
Returning to equation (\ref{GenLinCase}) this means that given $\theta$, there exists a distinguished {\it fundamental system}
\be{FUSYS}\wcirc_{1}(z|\theta),\ldots, \wcirc_{n}(z|\theta)\ee
such that
\be{asympwcirc}\wcirc_{j}\!\!(z|\theta)=e^{\Pi_j(z^{1/r})}z^{\rho_j} \Omega_j(z|\theta)\quad{\rm on~}S_\theta;\ee
$\Pi_j$ is a polynomial in $z^{1/r}$, $\rho_j$ is some complex number, and
$\Omega_j$ is a polynomial in $\log z$ with coefficients that have
asymptotic expansions in $z^{1/r}$ on $S_\theta$.
The triples $(\Pi_j,\rho_j,\Omega_j)$ are mutually distinct.
Again this system may vary from sector to sector, but only in the coefficients of $\Omega_j(z|\theta)$!
The leading terms of the $\Pi_j$ can be found among those of the $\int y_\nu\,dz.$

\subsection{\sc The Phragm\'en-Lindel\"of indicator}\label{PLI}
The following facts can be found in Lewin's/Levin's monographs~\cite{lewin,levin}. Let $f$ be an entire function of positive finite order $\varrho$ such that $\log M(r,f)=O(r^\varrho)$ as $r\to\infty$. Then
\be{PLIdef}h(\vartheta)=\limsup_{r\to\infty}\frac{\log|f(re^{i\vartheta})|}{r^{\varrho}}\quad(\vartheta\in\R)\ee
is called {\it Phragm\'en-Lindel\"of indicator} or just indicator of $f$ {\it of order} $\varrho$; it is continuous and always assumed to be extended to the real axis as a $2\pi$-periodic function.
Having $h(\vartheta)$ at hand, the Nevanlinna functions may be computed explicitly (for definitions and results we refer to Hayman~\cite{hayman} and Nevanlinna~\cite{nevanlinna}):
$$T(r,f)=m(r,f)\sim \frac{r^\varrho}{2\pi}\int_0^{2\pi}h^+(\vartheta)\,d\vartheta,[\footnote{$\psi_1(r)\sim \psi_2(r)$ as $r\to\infty$ means $\psi_1(r)/\psi_2(r)\to 1$.}]$$
$$m(r,1/f)\sim\frac{r^\varrho}{2\pi}\int_0^{2\pi}h^-(\vartheta)\,d\vartheta,\und N(r,1/f)\sim\frac{r^\varrho}{2\pi}\int_0^{2\pi}h(\vartheta)\,d\vartheta,$$
where as usual $x^+=\max\{x,0\}$ and $x^-=\max\{-x,0\}.$
In connection with our distinguished fundamental system (\ref{FUSYS}) only local indicators
\be{INDIloc}\ts h_j(\vartheta)=-\frac1\varrho\cos\big(\varrho\vartheta+2j\pi/(n-q)\big)~(0\le j<n-q)\und h_{n-q}(\vartheta)=0\ee
have to be considered.

\section{\sc Results on normal and sub-normal solutions}
The Airy integral has an asymptotic representation
$${\rm Ai}(z)\sim \frac1{2\sqrt\pi}e^{-\frac23 z^{3/2}}z^{-1/4}\sum_{j=0}^\infty c_jz^{-j/2}\quad{\rm even~on~}|\arg z|<\pi,$$
hence
\begin{itemize}
\item[1.] satisfies $\ds\frac{\log {\rm Ai}(z)}{z^{3/2}}\sim -\frac1\varrho$ as $z\to\infty$ on $|\arg z|\le\pi-\epsilon$ for every $\epsilon>0$;
\item[2.] has Phragm\'en-Lindel\"of indicator $h(\vartheta)=-\frac23\cos\big(\frac32\vartheta\big)$ on $|\vartheta|<\pi$;
\item[3.] has infinitely many zeros, all on the negative real axis.
\item[4.] $T(r,{\rm Ai})\sim\frac{8}{9\pi}r^{3/2}$ and $N(r,1/{\rm Ai})\sim\frac{4}{9\pi}r^{3/2}.$\end{itemize}

In the general case one cannot expect results of this high precision, but the following theorem seems to be a good approximation.

\subsection{\sc Theorem}\label{MAIN} {\it Any Laplace contour integral $\Lambda_L$  of order $\varrho=1+1/(n-q)\le 3/2$ {\rm (hence $n-q\ge 2$)}
\begin{itemize}
\item[1.] satisfies $\ds\frac{\log\Lambda_L(z)}{z^\varrho}\sim-\frac1\varrho$ on $|\arg z|\le \pi-\epsilon$  for every $\epsilon>0;$
\item[2.] has Phragm\'en-Lindel\"of indicator $h(\vartheta)=-\frac1\varrho\cos\big(\varrho\vartheta\big)$ on $|\vartheta|<\pi$;
\item[3.] has only finitely many zeros on $|\arg z-\pi|<\epsilon$ for every $\epsilon>0$;
\item[4.] $T(r,\Lambda_L)\sim\frac{1}{\pi\varrho^2}(1+|\sin(\pi\varrho)|)r^\varrho$ and $N(r,1/\Lambda_L)\sim\frac{1}{\pi\varrho^2}r^\varrho.$
\end{itemize}}\medskip

$\Lambda_L$ has `many' zeros by the last property; they are distributed over arbitrarily small sectors about the negative real axis.
The analogue of Theorem~\ref{MAIN} for  $q=n-1$ reads as follows.\medskip

{\sc \ref{MAIN}a.}  {\it Suppose $q=n-1$. Then either $\Lambda_L(z)=e^{-z^2/2}P(z)$, $P$ some non-trivial polynomial, or else
the following is true:
\begin{itemize}
\item[1.] $\ds\frac{\log\Lambda_L(z)}{z^2}\sim -\frac 12$ on $|\arg z|<\frac34\pi-\epsilon$;
\item[2.] $h(\vartheta)=-\frac12\cos(2\vartheta)$ on $|\vartheta|<\frac34\pi$ and $h(\vartheta)=0$  on $\frac34\pi\le|\vartheta|\le \pi$;
\item[3.] $\Lambda_L$ has only finitely many zeros on $|\arg z|\le\frac34\pi-\epsilon$, and at most\\ $o(r^2)$ zeros on $|\arg z-\pi|<\frac\pi 4-\epsilon$, $|z|<r$, for every $\epsilon>0;$
\item[4.] $T(r,\Lambda_L)\sim\frac{1}{2\pi}r^2$ and $N(r,1/\Lambda_L)\sim\frac{1}{4\pi}r^2$.
\end{itemize}}\medskip

In the second case, $\Lambda_L$ again has `many' zeros; they are distributed over arbitrarily small sectors about $\arg z=\pm\frac34\pi$.\medskip

The proof of Theorem~\ref{MAIN} and \ref{MAIN}{\sc a} will be given in section~\ref{Beweis}.

\subsection{\sc Subnormal solutions}\label{sns} Solutions of maximal order $\varrho=1+1/(n-q)$ always exist, while the existence of so-called {\it sub-normal} solutions
having order  $\varrho=1$,  $\varrho=1-1/p$, and  $\varrho=0$ (polynomials), respectively, remains open; necessary but by far not sufficient conditions are   $q>p$,  $p>1$, and $p=1$,
respectively. Sufficient condition are coupled with the poles of $Q_0/Q_1$ and their residues. We note that at $t=0$, $Q_0/Q_1$ either has a pole of order $p\ge 1$ or is regular ($p=0)$.
We have to distinguish two cases.\medskip

{\sc \ref{sns}a~Theorem.} {\it Let $t=0$ be a pole of $Q_0/Q_1$ of order $p\ge 1$ with integer residue $\lambda_0$. Then the solution
$$w(z)=\Res_{0}\big[\phi_L(t)e^{-zt}\big]$$
has order of growth $\varrho=1-1/p$ if $t=0$ is an essential singularity of $\phi_L$, that is, if $p>1$, and otherwise $(p=1)$ is either a non-trivial polynomial of degree $\lambda_0\ge 0$ or else vanishes identically $(\lambda_0<0)$.}\medskip

{\sc \ref{sns}b~Theorem.} {\it Let $t=t_0\ne 0$ be a pole of $Q_0/Q_1$ of order $m$ with integer residue $\lambda_0$.
Then the solution
$$w(z)=\Res_{t_0}\big[\phi_L(t)e^{-zt}\big]$$
has the form
$$w(z)=e^{-t_0z}W(z),$$
where $W$ is a transcendental entire function of order of growth at most $1-1/m$ if $m>1$, and otherwise is a polynomial of degree $\lambda_0\ge 0$ or vanishes identically $(\lambda_0<0)$.}\medskip

{\sc Proof} of \ref{sns}{\sc a}. First assume that $t=0$ is a simple pole of $Q_0/Q_1$. Then $p=1$ (since $Q_0(0)=a_0\ne 0$) and
$$\phi_L(t)e^{-zt}=\frac{H(t)e^{-zt}}{(t-t_0)^{1+\lambda_0}}$$
holds, where $H$ is holomorphic at $t=0$ and $H(0)\ne 0$. Then (assuming $\lambda_0\ge 0$)
$$w(z)=\sum_{k=0}^{\lambda_0}\frac{H^{(\lambda_0-k)}(0)}{k!(\lambda_0-k)!}(-z)^k$$
is a polynomial of degree $\lambda_0$. If, however $t=0$ is essential for $\phi_L$, then $p>1$ and
$$|\phi_L(t)|\le \frac{Ae^{B|t|^{-p+1}}}{|t|^{1+\lambda_0}}$$
holds on $|t|<\delta,$ hence
$$|w(z)|\le \frac1{2\pi }\int_{|t|=|z|^{-1/p}}\frac{Ae^{B|t|^{-p+1}+|z||t|}}{|t|^{1+\lambda_0}}\,|dt|\le A|z|^{\lambda_0/p}e^{(B+1)|z|^{1-1/p}}$$
holds as $z\to\infty$. This yields $\varrho(w)\le 1-1/p$. Since $w$ is not a polynomial, the only possibility is $\varrho=1-1/p$. This proves part {\sc a}.\medskip

To deal with \ref{sns}{\sc b} write
$$w(z)=\Res_{t_0}\big[\phi_L(t)e^{-zt}\big]=e^{-t_0z}\Res_{t_0}\big[\phi_L(t)e^{-z(t-t_0)}\big]=e^{-t_0z}W(z).$$
Then by part {\sc a} of the proof, either $W$ is a polynomial of degree $\lambda_0\ge 0$ or vanishes identically $(\lambda_0<0$) or is a transcendental entire function of order of growth at most
$1-1/m$; the value zero is a Borel or even Picard exceptional value of $w$.\hfill\Ende\medskip

{\sc Remark.} By $w=e^{-t_0 z}W$ ($t_0\ne 0$), equation (\ref{GenLinCase}) is transformed into
\be{WDGL}W^{(n)}+\sum_{j=0}^{n-1}(A_j+B_jz)W^{(j)}=0\ee
with $A_\nu+B_\nu z=\sum_{j=\nu}^n{j\choose \nu}(-t_0)^{j-\nu}(a_j+b_jz);$
in particular, $B_q=b_q$ and $A_0+B_0z=Q_0(t_0)+Q_1(t_0)z$
hold, hence $q$, $b_q$, and the maximal order of growth $1+1/(n-q)$ are invariant under this transformation, but not necessarily the index $p$. If in Theorem \ref{sns}{\sc b},
$t_0$ is not a zero of $Q_0(t)$ and $m>1$, then in (\ref{WDGL}) we have $B_m\ne 0$, but $B_j=0$ for $j<m$, hence $W$ has the order $1-1/m$ by part {\sc a}.

\section{\sc Rotational symmetries}\label{RotSymm}

The functions ${\rm Ai}(ze^{-2\nu\pi i/3})$ also solve Airy's equation and coincide up to some non-zero factor with the Laplace contour integrals
$$\frac1{2\pi i}\int_{\CC_{\frac{2\nu+1}3\pi,\frac{2\nu-1}3\pi}}\!\!\!\!\!\!\!\!\!\!\!\! e^{t^3/3}e^{-zt}\,dt;$$
obviously the sum over the three integrals vanishes identically, while any two of these functions are linearly independent; this is well-known, and will be proved below in a more general context. In the general case there are two obvious obstacles:\medskip

{\it \qquad 1.\quad $\Lambda_L(z e^{-2\nu\pi i/(n-q+1)})$ is not necessarily a solution to $L[w]=0$.}

Nevertheless one may consider the contour integral solutions
\be{varpidef}\Lambda_{\nu}(z)=\frac1{2\pi i}\int_{\CC_{R,\theta_{2\nu+1},\theta_{2\nu-1}}}\!\!\!\!\!\!\!\!\!\!\!\!\phi_L(t)e^{-zt}\,dt\quad(\ts\theta_k=\frac{k\pi}{n-q+1})\ee
and their Phragm\'en-Lindel\"of indicators $h_\nu$ for $0\le \nu\le n-q$, but has to take care since $\phi_L$ may be many-valued.
The next theorem is just a reformulation of Theorem~\ref{MAIN} and \ref{MAIN}{\sc a}. We assume the indicator $h_0$ of $\Lambda_0=\Lambda_L$ be extended to the real line as a $2\pi$-periodic function.
Although $\Lambda_\nu(z)$ is different from $\Lambda_L(ze^{-2\nu\pi i/(n-q)})$ in general, both functions share their main properties as are stated in Theorem~\ref{MAIN} and \ref{MAIN}{\sc a} for $\Lambda_L$.

\subsection{\sc Theorem} {\it Each contour integral $(\ref{varpidef})$ of order $\varrho=1+1/(n-q)$ has the same properties{\rm [\footnote{concerning asymptotics, Phragm\'en-Lindel\"of indicator, distribution of zeros, and Nevanlinna functions}]} as the corresponding function $\Lambda_L(e^{-2\nu\pi i/(n-q+1)}z)$ in {\rm Theorem}~$\ref{MAIN}$ and $\ref{MAIN}${\sc a}.}

\medskip If $\phi_L$ has no singularities, the sum
\be{varpisum}\Lambda(z)=\Lambda_L(z)+\Lambda_{1}(z)+\cdots + \Lambda_{n-q}(z)=-\frac1{2\pi i}\int_{|t|=R}\phi_L(t)e^{-zt}\,dt\ee
vanishes identically. On the other hand, poles $t_k$ of $Q_0/Q_1$ are poles or essential singularities of $\phi_L$; if these singularities are non-critical, that is, if the residues $\lambda_k=\Res_{t_k}[Q_0/Q_1]$ are integers, the following holds.

\subsection{\sc Theorem} {\it Suppose the residues $\lambda_k$ are integers. Then the sum $\Lambda$ either vanishes identically or may be written as a linear combination of sub-normal solutions $$\Res_{t_k}[\phi_L(t)e^{-zt}]=e^{-t_k z}W_k;$$
$W_k$ has order of growth less than one.}

{\it\medskip\qquad 2.\quad $\phi_L$ may be many-valued on  $\C\setminus\{$ {\rm poles of} $ Q_0/Q_1\}$.}

This is the case if some of the residues $\lambda_k$ are not integers.
Nevertheless $\phi_L$ may be single-valued on $|t|>R_0$: if $t$ goes around once the positively oriented circle $|t|=R>R_0$, $(t-t_k)^{\lambda_k}$
takes the value $(t-t_k)^{\lambda_k}e^{2\pi i\lambda_k}$, hence analytic continuation of $\phi_L(t)$ along $|t|=R$ yields the value $\phi_L(t)e^{-2\pi i\sum_{k}\lambda_k}$ (see the general form
(\ref{Generalphi}) of $\phi_L$), and $\phi_L$ is single-valued on $|t|>R_0$ if $\sum_{k}\lambda_k$ is an integer. Assuming this,
$\Lambda$ is a sub-normal solution of order at most one. This follows from  (\ref{varpisum}) and $|\phi_L(t)e^{-zt}|\le Ae^{R|z|}$ on $|t|=R,$ hence
$$|\Lambda(z)|\le ARe^{R|z|}.$$
We have thus proved

\subsection{\sc Theorem} {\it If the sum of residues $\sum_k\lambda_k$ is an integer, the sum $\Lambda$ either vanishes identically or is a sub-normal solution to $L[w]=0$.}\medskip

The next theorem is concerned with the linear space spanned by the functions $\Lambda_\nu$, $0\le\nu\le n-q$.

\subsection{\sc Theorem} {\it Any set of $n-q$ functions $\Lambda_\nu$ is linearly independent.}\medskip

{\sc Proof.} Set
$$y(\vartheta)=\max_{0\le\nu\le n-q} h_\nu(\vartheta)=\max_{0\le\nu\le n-q} h_0(\vartheta-\theta_{2\nu})$$
and $I_\nu=(-\pi+\theta_{2\nu},-\pi+\theta_{2\nu+2})$, $0\le \nu\le n-p$. Then
$$h_k(\vartheta)=h_{k+1}(\vartheta)=h(\vartheta)>h_j(\vartheta)\quad{\rm holds~on~}I_{n-p-k}$$
for every $j\ne k,k+1$ mod $(n-p+1)$. Now let
$$\sum_{\nu=0}^{n-p}c_\nu \Lambda_\nu$$
be any non-trivial linear combination of the trivial solution, if any, and assume $c_k\ne 0$ for some $k$. Then
$$h_k(\vartheta)\le\max\{h_\nu(\vartheta): c_\nu\ne 0,~\nu\ne k\}$$
holds on the other hand, and
$$h_k(\vartheta)=h(\vartheta)$$
on $I_{n-p-k}\cup I_{n-p-k+1}$ on the other.
This implies $c_{k\pm 1}\ne 0$, hence $c_\nu\ne 0$ for every $\nu$, and proves that any collection of $n-q$ functions $\Lambda_\nu$ is linearly independent. \hfill\Ende\medskip

\begin{center}\includegraphics[scale=0.3]{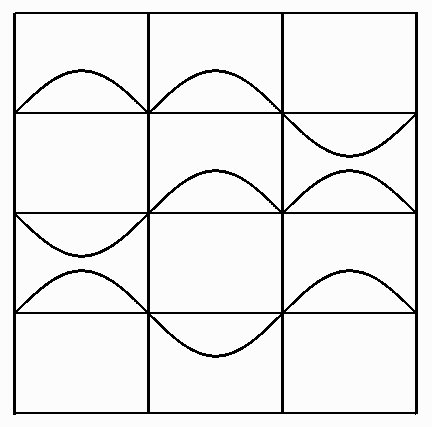}\quad\includegraphics[scale=0.3]{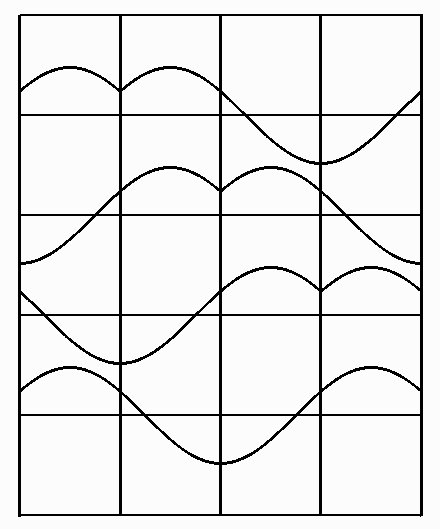}\quad\includegraphics[scale=0.3]{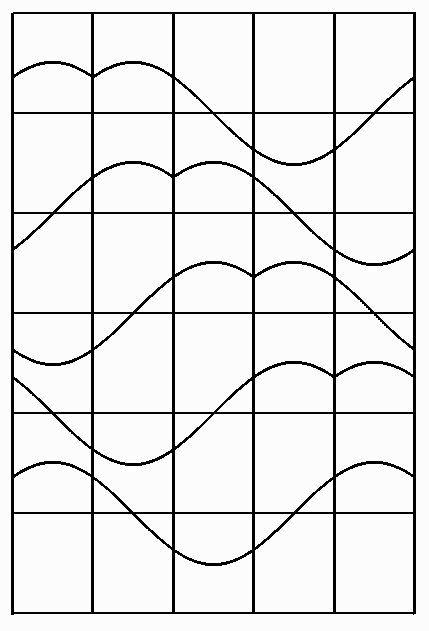}

\small{\sc Figure.} The indicators $h_k$ on $[-\pi,\pi]$ ($0\le k\le n-q$, from bottom to top), and the intervals $I_\nu$ ($0\le \nu\le n-q$, from left to right) for $n-q=2, 3, 4$.\end{center}

\medskip{\sc Remark.} Note that $\phi_L$ has no singularities at all only in case of
$$L[w]=w^{(n)}+\sum_{j=0}^{n-1}a_jw^{(j)}+(-1)^{n+1}zw=0$$
(with our normalisation). The most simple three-term example is equation~(\ref{FAM1}). The Laplace solutions to
$$w^{(n)}+(-1)^{n+1}zw=0$$
seem to be qualified candidates to be named special functions.
We do not know what happens if  $\sum_k\lambda_k$ is not an integer. Can one guarantee the existence of sub-normal solutions?
Is the sum $\Lambda$ itself sub-normal?

\section{\sc More Examples}

\begin{itemize}
\item[1)] The differential equation $w^{(5)}-zw'''+7w''+zw'-2w=0$\\
\noindent ($n=5,~q=3,~p=1$, $Q_0(t)=-t^5+7t^2-2$, $Q_1(t)=t^3-t$) has the fundamental system
$$\begin{aligned}
\qquad w_1=&~\Lambda_L(z)=\frac1{2\pi i}\int_{\CC_{\frac\pi 3}}\frac{e^{t^3/3+t-zt}}{t^3(t-1)^3(t+1)^4}\,dt\cr
w_2=&\Lambda_1 {\rm ~or~} \Lambda_2\cr
w_3=&~(2z^3-39z^2+264z-635)e^z\cr
w_4=&~(z^2+6z+15)e^{-z}\cr
w_5=&~z^2+7.\end{aligned}$$
The set of functions $\Lambda_0$, $\Lambda_1$, and $\Lambda_2$ is linearly dependent.

\item[2)] The differential equation $w^{(6)}-zw^{(4)}-zw''+w=0$\\
\noindent with $n=6$, $q=4$, $p=2$, $Q_0(t)=t^6+1$ and $Q_1(t)=-t^4-t^2$, hence
$$\phi_L(t)=\frac1{t^4+t^2}\exp\Big[\frac{t^3}3-t-\frac 1t\Big],$$
has solutions $w_{1,2}=e^{\pm iz}$ (associated with the poles at $t=\mp i$) of order $\varrho=1$ and three solutions of order $\varrho=3/2$, namely $\Lambda_0$, $\Lambda_1$, and $\Lambda_1$.
The residue theorem yields one more solution $w_3(z)=\Res_0\big[\phi_L(t)e^{-zt}\big]$ that has order of growth $2/3$.
It is not clear how complete the five solutions just discussed to obtain a basis since
$$\Lambda_0(z)+\Lambda_1(z)+\Lambda_2(z)=\frac i2e^{i(z+1/3)}-\frac i2e^{-i(z+1/3)}-w_3(z)$$
is a linear combination of $w_1,$ $w_2$, and $w_3$.
\medskip

\item[3)] Section 6 in \cite{GHW} is devoted to the differential equation
$$v^{(4)}-zv'''-v=0$$
(unfortunately not with our normalisation). The authors found three contour integral solutions
denoted $H(z)$, $H(-z)$, and $U(z)=H(z)+H(-z)$ of order $3/2$; any two of them are linearly independent.
The transformation $w(z)=v(iz)$ transforms this equation into
$$w^{(4)}+zw'''-w$$
with $Q_0(t)=t^4-1$, $Q_1(t)=-t^3$, and $\phi_L(t)=\frac1{t^3}\exp\Big[\frac{t^2}2+\frac1{2t^2}\Big]$.
The solutions
$$\Lambda_j(z)=\frac1{2\pi i}\int_{\CC_{1,\frac{2j-1}2\pi,\frac{2j+1}2\pi}}\!\!\!\!\!\!\!\!\!\!\!\!\phi_L(t)e^{-zt}\,dt\quad(j=0,1)$$
have order of growth $2$, and the sum
$$-\Lambda_0(z)-\Lambda_1(z)=\frac1{2\pi i}\int_{|t|=1}\phi_L(t)e^{-zt}\,dt=\Res_0[\phi_L(t)e^{-zt}]$$
is sub-normal with order of growth $2/3$ ($n=4, q=p=3$); it corresponds to the solution $U$ in \cite{GHW} and has power series expansion
$$\sum_{m=0}^\infty \frac{C_m}{2^{m-1}m!}z^{2m}\quad{\rm with}\quad C_m=\sum_k\frac{1}{4^k k!(k+m-1)!};$$
$k$ runs over every non-negative integer such that $k+m-1\ge 0$.
\medskip

\item[4)] To every $n\ge 2$ there exist a unique family of differential equations
$$w^{(n)}+zw^{(n-1)}+\ds\sum_{j=0}^{n-2}(a_j+b_jz)w^{(j)}=0$$
depending on $n-2$ parameters with solution $w=e^{-z^2/2}$.

\item[5)] The differential equation $w^{(4)}+(z-1)w'''-8w''-zw'+2w=0$\\
\noindent ($n=4, q=3, Q_0(t)=t^4+t^3-8t^2+2, Q_1(t)=-t^3+t$) has the fundamental system
$$\begin{aligned}
w_1(z)=&~\Lambda_L(z)=\frac1{2\pi i}\int_{\CC_{2,\frac\pi 2}}\frac{e^{t^2/2+t-zt}}{t^3(t-1)^3(t+1)^4}\,dt\cr
w_2(z)=&~(2z^3-27z^2+150z-324)e^z\cr
w_3(z)=&~(z^2+6z+14)e^{-z}\cr
w_4(z)=&~z^2+8.
\end{aligned}$$
$\Lambda_0+\Lambda_1$ is a non-trivial linear combination of $w_2, w_3,$ and $w_4$.\medskip

\item[6)] The differential equation $w^{(5)}-(z+1)w'''+w''+(z+1-\lambda-2\mu)w'+\lambda w=0$
with
$Q_0(t)=-t^5+t^3+t^2+(\lambda+2\mu-1)t+\lambda$, $Q_1(t)=t^3-t$ ($\lambda$ and $\mu$ arbitrary)
has the distinguished solutions $\Lambda_0=\Lambda_L$, $\Lambda_1$ and $\Lambda_2$. Although
$$\phi_L(t)=\frac{e^{t^3/3}}{t^{1-\lambda}(t-1)^{1+\lambda+\mu}(t+1)^{2-\mu}}$$
may have transcendental singularities at $t=0,1,-1$, $\phi_L$ is single-valued on $|t|>1$, and $\Lambda_0(z)+\Lambda_1(z)+\Lambda_2(z)$
has order of growth at most $1$.
\end{itemize}

\section{Preparing the proof of theorem~\ref{MAIN}}

Our proof of Theorem~\ref{MAIN} will be based on the following

\subsection{\sc Proposition}\label{GHWPropo} {\it Suppose
\be{Localphi}\phi(t)=\exp\Big[\frac{t^{m+1}}{m+1}+\psi(t)\Big]\ee
is holomorphic on $\{t:|t|>R_0,~|\arg t|<\pi\}$, where $\psi$ satisfies $|\psi(t)|\le C|t|^{m}$ and $m\ge 1$ is an integer. Then
the Phragm\'en-Lindel\"of indicator of
$$f(z)=\frac1{2\pi i}\int_{\CC_{R,\frac\pi{m+1}}}\phi(t)e^{-zt}\,dt$$
of order $\varrho=1+1/m$ is negative on
\be{MaxSektor}|\theta|<\frac{m\pi}{2(m+1)}.\ee}

{\sc Remark.} Actually Proposition~\ref{GHWPropo} was stated and proved in \cite{GHW} (Theorem~3) for
$$\phi(t)=\exp\Big[\frac{t^{m+1}}{m+1}+\frac{(-1)^{k+1}b t^{k+1}}{k+1}\Big]\quad(0<k<m),$$
without reference to the Phragm\'en-Lindel\"of indicator. Examining the proof shows that it works for any $\phi$ given by (\ref{Localphi}) and  $m\ge 2$, but {\it not} for $m=1$,
which was out of sight in \cite{GHW}. For the proof of Theorem~\ref{MAIN} we only need $m\ge 2$ ($m=n-q$ in our notation), but for the addendum to this theorem we need the case $m=1$.
\medskip

{\sc Proof} of Proposition~\ref{GHWPropo} for $m=1$. Our object now is
\be{OBJEKT}\frac1{2\pi i}\int_{\CC_{R,\frac\pi{2}}}\phi(t)e^{-zt}\,dt\ee
with
$$\phi(t)=\exp\Big[\frac{t^2}2+\psi(t)\Big]\und|\psi(t)|\le C|t|.$$
First of all we notice that $R>R_0$ may take any value since for $R_1\ge R_2>R_0$, the simple closed curve $\CC_{R_2,\frac\pi{2}}\circleddash\CC_{R_1,\frac\pi{2}}$ is contained in the domain of $\phi$.
We choose $R=\epsilon |z|$, where $\epsilon>0$ depending on $\theta=\arg z$ will be determined during the proof. Secondly we choose $\lambda>1$ such that $\lambda\theta<\pi/4$
and show that for $0<\theta<\pi/4$, say, the contour $\CC_{\epsilon|z|,\frac\pi 2}$ may be replaced with $\CC_{\epsilon|z|,\frac\pi 2-\lambda\theta,-\frac\pi 2}$
(of course, nothing has to be done if $\theta=0$). To this end we have to show that the integral over the arc $\Gamma_r: t=ire^{-i\vartheta}$, $0\le\vartheta\le\lambda\theta$
vanishes in the limit $r\to\infty$. This follows from
$$\begin{aligned}
\Big|\frac1{2\pi i}\int_{\Gamma_r}\phi(t)e^{-zt}\,dt\Big|\le &~\frac r{2\pi}\int_0^{\lambda\theta}\exp\Big[-\frac{r^2}2\cos (2\vartheta)+|z|r+Cr\Big]\,d\vartheta\cr
\le&~\frac r8\exp[-\frac{r^2}2\cos (2\lambda\theta)+r(|z|+C)\Big]\to 0~{\rm as~} r\to\infty\end{aligned}$$
since $\cos(2\lambda\theta)>0$. To finish the proof we have to estimate the integrals
$$\int_{\CC}|\phi(t)e^{-zt}|\,|dt|$$
over
\begin{itemize}\item[1.]  $\CC: t=-i\epsilon|z|\tau$, $1\le\tau<\infty$;
\item[2.]$\CC: t=i\epsilon|z|\tau e^{-i \lambda\theta}$, $1\le\tau<\infty$;
\item[3.]$\CC: t=\epsilon|z|e^{i\vartheta}$, $-\pi/2\le\vartheta\le \pi/2-\lambda\theta$.
\end{itemize}\medskip

\begin{itemize}\item[ad 1.] From $\Re t^2=-\epsilon^2|z|^2\tau^2\le -\epsilon^2|z|^2\tau$, $\Re(-zt)=-\epsilon|z|^2\tau\sin\theta\le 0$ (since $0\le\theta<\pi/4$), and $|dt|=\epsilon|z|\,d\tau$ we obtain the upper bound
$$\qquad\epsilon|z|\int_1^\infty\exp\Big[-\frac12\epsilon|z|(\epsilon|z|-2C)\tau\Big]\,d\tau\le\frac1Ce^{-\epsilon^2|z|^2/4}\quad(|z|\ge 4C/\epsilon).$$
Note that every $\epsilon>0$ works and $\lambda$ doesn't appear.\medskip

\item[ad 2.] Here it follows from $\Re t^2=-\epsilon^2|z|^2\tau^2\cos(2\lambda\theta)\le -\epsilon^2|z|^2\tau\cos(2\lambda\theta)$ and
$\Re(-zt)=-|z|^2\tau^2\sin((\lambda-1)\theta)\le 0$ that
$$\epsilon|z|\int_1^\infty\exp\Big[-\frac12\epsilon|z|(\epsilon|z|\cos(2\lambda\theta)-2C)\tau\Big]\,d\tau \le\frac 1Ce^{-\epsilon^2|z|^2/4}$$
is an upper bound provided $\ds|z|\ge \frac{4C}{\epsilon\cos(2\lambda\theta)}$. Again we note that every $\epsilon>0$ works and $\lambda$ does not play an essential role.\medskip
\item[ad 3.] We first assume $0<\theta<\lambda\theta<\pi/4$. From  $\Re t^2\le\epsilon^2|z|^2$, $\Re (-zt)=-\epsilon|z|^2\cos (\vartheta+\theta)$, and $|dt|=\epsilon|z|\,d\vartheta$
we get the bound
$$\epsilon|z|\int_{-\pi/2}^{\pi/2-\lambda\theta}\exp\Big[\frac12\epsilon|z|^2(\epsilon-2\cos(\vartheta+\theta))+\epsilon C|z|\Big]\,d\vartheta$$
To control the term $\cos(\vartheta+\theta)$ note that
$$\qquad-\pi/2<-\pi/2+\theta\le \vartheta+\theta\le \pi/2-\lambda\theta+\theta=\pi/2-(\lambda-1)\theta<\pi/2,$$
hence
$\cos(\vartheta+\theta)\ge\min\{\sin\theta, \sin((\lambda-1)\theta)\}=\kappa(\theta)>0;$
here we need $\lambda$! This yields the upper bound
$$\pi\epsilon|z|\exp\Big[\frac12\epsilon|z|^2(\epsilon-2\kappa(\theta))+\epsilon C|z|\Big]$$
for our integral, and all we have to do is to choose $\epsilon=\kappa(\theta)$ to obtain the bound $\pi\epsilon|z|e^{-\epsilon^2|z|^2/4}$ for $|z|\ge 4C/\epsilon$.
It remains to discuss the case $\theta=0$, hence $z=x>0$, where we can work with $\epsilon=1$. From
$$\qquad\qquad\Re\Big[\frac{t^2}2-xt\Big]=\frac12\Re(t-x)^2-\frac12x^2=2x^2(\cos\vartheta-1)\cos\vartheta-\frac12x^2\le-\frac12x^2$$
we obtain the very last upper bound
$$\pi\exp\Big[-\frac12x^2+Cx\Big]\le \pi\exp\Big[-\frac14x^2\Big]\quad(x\ge 4C).$$
\end{itemize}
This proves Proposition~\ref{GHWPropo} for $m=1$. To verify the case $m\ge 2$ the reader is referred to \cite{GHW}.\hfill\Ende
\medskip

\section{\sc Proof of Theorem~\ref{MAIN} and \ref{MAIN}{\sc a}}\label{Beweis}
Let $h$ denote the Phragm\'en-Lindel\"of indicator of our special solution $\Lambda_L$. To prove
\be{hequalsh0}h(\vartheta)=-\frac1\varrho\cos(\varrho\vartheta)\ee
on $(-\frac\pi{2\varrho}\le\vartheta\le \frac\pi{2\varrho})$ we use the following facts
taken from \cite{levin}, p.\ 53 ff., which hold for arbitrary indicators of order $\varrho$.

\begin{itemize}
\item[1.] $h$ has one-sided derivatives everywhere, and $h'(\vartheta-)\le h'(\vartheta+)$ holds;
\item[2.]  $h$ is $\varrho$-{\it trigonometrically convex};
\item[3.] $h(\varphi)+h(\varphi+\pi/\varrho)\ge 0$ holds for every $\vartheta$.\end{itemize}
Since, however, $h(\pm\frac\pi{2\varrho})\le 0$, the third property implies $h(\pm \frac\pi{2\varrho})=0$, and the second leads to
$$h(\vartheta)=h(0)\cos(\varrho\vartheta)\quad{\rm on~}\ts |\vartheta|\le \frac\pi{2\varrho}$$
(again \cite{levin}, p.\ 53 ff.). Then $h(0)=-1/\varrho$ follows from the fact that the only local indicators available are
given by (\ref{INDIloc}).
By the third property (note $h'(\frac\pi{2\varrho})=1$), (\ref{hequalsh0}) remains true on $(\frac\pi{2\varrho},\theta]$, where
$\theta>\frac\pi{2\varrho}$ is the smallest number such that
$$h_0(\theta)=h_k(\theta)\und h'_0(\theta)<0 {\rm~and~}h'_k(\theta)>0$$
holds for some $k\ne 0$. This happens soonest at $\theta=\pi$. 
Finally, the claim about the distribution of zeros follows from the fact that
$$h(\vartheta)=\lim_{r\to\infty}\frac{\log|\Lambda_L(re^{i\vartheta})|}{r^\varrho}$$
holds uniformly on $|\vartheta|\le\pi-\epsilon$, and
$$\frac{\log\Lambda_L(z)}{z^\varrho}\to-\frac1\varrho\quad{\rm as}~z\to\infty~{\rm on~}|\arg z|\le\pi-\epsilon.$$
Thus all but finitely many zeros are contained in $|\arg z-\pi|<\epsilon$ for every $\epsilon>0$.
This proves Theorem~\ref{MAIN}, and the arguments  may be repeated step-by-step to prove Theorem~\ref{MAIN}{\sc a} until $\vartheta=\pm\frac34\pi$
is reached. Then either $h(\vartheta)=-\frac12\cos(2\vartheta)$ holds on $[-\pi,\pi]$ or only on $[-\frac34\pi,\frac34\pi]$, while $h(\vartheta)=0$ on $\frac34\pi\le|\vartheta|\le\pi$.
Then
$$\ds\frac{\log\Lambda_L(z)}{z^2}\to-\frac12\quad{\rm as~}z\to\infty$$
either holds on the whole plane or on $|\arg z|\le\frac34\pi-\epsilon$. In the first case $\Lambda_L$ has only finitely many zeros and $\Lambda_L(z)=e^{-z^2/2}P(z)$ holds.
In the second case, $\Lambda_L$ has only finitely many zeros on $|\arg z|\le\frac34\pi-\epsilon$, and from
$$\log|\Lambda_L(re^{i\vartheta})|=o(r^2)$$
as $r\to\infty$, uniformly on $|\vartheta|\le\frac\pi4-\epsilon$, it follows that $\Lambda_L$ has $o(r^2)$ (probably only $O(r)$) zeros on $|\arg z-\pi|\le\frac\pi4-\epsilon$, $|z|\le r$  (cf.\ \cite{lewin}, p.\ 150). The asymptotic formulae for $T(r,\Lambda_L)$ and $N(r,1/\Lambda_L)$ are obvious.\hfill\Ende

\bigskip{\footnotesize Fakult\"at f\"ur Mathematik, Technische Universit\"at Dortmund\\
{\sf stein@math.tu-dortmund.de, http://www.mathematik.tu-dortmund.de/steinmetz/}\\
postal address: Beethovenstrasse 17, D-67360 Lingenfeld, Germany}
\end{document}